\documentclass[titlepage,12pt]{amsart}
 \usepackage{amsfonts}
 \usepackage{amssymb}
 \usepackage{amsmath,amsxtra,amsthm}

\usepackage[usenames]{color}

\usepackage[dvips]{graphicx}

\usepackage[all,matrix,arrow,curve]{xy}

\usepackage{mathrsfs}

\usepackage{enumitem}

\DeclareFontFamily{OT1}{pzc}{}
\DeclareFontShape{OT1}{pzc}{m}{it}{<-> s * [1.15] pzcmi7t}{}
\DeclareMathAlphabet{\mathpzc}{OT1}{pzc}{m}{it}

\DeclareSymbolFont{bbold}{U}{bbold}{m}{n}
\DeclareSymbolFontAlphabet{\mathbbold}{bbold}

\newcommand\iso{\xrightarrow{
   \,\smash{\raisebox{-0.45ex}{\ensuremath{\sim}}}\,}}

\newtheorem{theorem}{Theorem}
\newtheorem{deff}[theorem]{Definition}

\newtheorem{proposition}[theorem]{Proposition}
\newtheorem{example}[theorem]{Example}
\newtheorem{lemma}[theorem]{Lemma}

\newtheorem{cor}[theorem]{Corollary}

\newtheorem{rem}[theorem]{Remark}

\newtheorem*{theorem*}{Theorem}

\numberwithin{equation}{section}

\renewcommand{\proof}{{\bf Proof.}~}

\newcommand{\mto}{\mapsto}
\newcommand{\bqa}{\begin{eqnarray}}
\newcommand\eqa {\end{eqnarray}}
\newcommand{\beq}{\begin{eqnarray}}
\newcommand{\beqn}{\begin{eqnarray}\nonumber}
\newcommand{\eeq}{\end{eqnarray}}
\newcommand{\be}{\begin{array}}
\newcommand{\ee}{\end{array}}

   \newcommand\vf\varphi

 \newcommand{\End}{\mathrm{End}}
 \newcommand{\Ker}{\mathrm{Ker}}
 \newcommand{\Id}{\mathrm{Id}}

 \newcommand{\md}{\mathrm{d}}

 \newcommand{\cB}{\mathcal{B}}

 \newcommand{\calG}{\mathcal{G}}

 \newcommand{\N}{{\mathbb N}}

\title[Riemann-Cartan groupoids and quadratic Lie algebroids]{Integration of quadratic Lie algebroids to Riemannian Cartan-Lie groupoids }%
\author{Alexei Kotov}
\address{Alexei Kotov: Faculty of Science, University of Hradec Kralove, Rokitanskeho 62, Hradec Kralove
50003, Czech Republic}
\email{oleksii.kotovAuhk.cz}

\author{Thomas Strobl}
\address{Thomas Strobl: UMI CNRS2924 Instituto de Matem\'atica Pura e Aplicada (IMPA), Estrada Dona Castorina 110, Rio de Janeiro, 22460-320, Brasil and \newline
Institut Camille Jordan, Universit\'e Claude Bernard Lyon 1, Universit\'e de Lyon,
43 boulevard du 11 novembre 1918, 69622 Villeurbanne Cedex, France}
\email{thomasATimpa.br and stroblATmath.univ-lyon1.fr}

\date{December 2017}

\begin{document}

\thispagestyle{empty}

  \begin{abstract}

Cartan-Lie algebroids, i.e.~Lie algebroids equipped with a compatible connection, permit the definition of an adjoint representation, on the fiber as well as on the tangent of the base. We call (positive) quadratic Lie algebroids, Cartan-Lie algebroids with ad-invariant (Riemannian) metrics on their fibers and base $\kappa$ and $g$, respectively. We determine the necessary and sufficient conditions for a positive quadratic Lie algebroid to integrate to a Riemmanian Cartan-Lie groupoid. Here we mean a Cartan-Lie groupoid $\calG$ equipped with a bi-invariant and inversion invariant metric $\eta$ on $T\calG$ such that it induces by submersion the metric $g$ on its base and its restriction to the $t$-fibers coincides with $\kappa$.

\vspace{20mm}
\noindent \keywordsname : Lie algebroids, Lie groupoids, Cartan connections, multiplicative distributions, jet spaces
and jet bundles, Riemannian submersions

\vspace{4mm}

\vspace{2.5cm}

\tableofcontents

  \end{abstract}

\setcounter{page}{0}

\maketitle

\clearpage

\setcounter{page}{1}

  \def\sp{\mathfrak sp}
  \def\sll{\mathfrak sl}
  \def\g{{\mathfrak g}}
  \def\gl{{\mathfrak gl}}
  \def\su{{\mathfrak su}}
  \def\so{{\mathfrak so}}
  \def\sll{{\mathfrak sl}}
  \def\h{{\mathfrak h}}

  \def\P{{\mathbb P}}
  \def\H{\mathbb H}

   \def\a{\alpha}
   \def\b{\beta}
   \def\t{\theta}
   \def\la{\lambda}
   \def\e{\epsilon}
   \def\ga{\gamma}
   \def\de{\delta}
   \def\De{\Delta}
   \def\om{\omega}
   \def\Om{\Omega}

   \def\i{\imath}

 \def\gr{\g^{\scriptscriptstyle\mathrm{gr}}}
 \def\godd{\g_{\scriptscriptstyle 1}}
 \def\geven{\g_{\scriptscriptstyle 0}}
 \def\grodd{\gr_{\scriptscriptstyle 1}}
 \def\greven{\gr_{\scriptscriptstyle 0}}

  \def\sst{\scriptscriptstyle}

  \def\sot{{\;{\scriptstyle \otimes}\;}}
  \def\st{{\sst\times}}

  \def\df{{\sst\mathrm{def}}}

\def\scrE{\mathscr{E}}
  \def\scrF{\mathscr{F}}
  \def\scrJ{\mathscr{J}}
  \def\scrR{\mathscr{R}}


\section*{Introduction}

\noindent In earlier papers \cite{Kotov-Strobl14,Kotov-Strobl16_Geom1} we found a compatibility condition for a Lie algebroid over a Riemannian base, when the bundle of the algebroid is equipped with a connection. In \cite{Kotov-Strobl14} we were interested in gauging symmetries of the standard sigma model whose fields take values in a Riemannian manifold; it turned out that the Lie algebra of Killing symmetries can be replaced by a more general Lie algebroid, which we called a Killing Lie algebroid.
Even if the target Riemannian manifold does not admit any Killing vector field, such an algebroid may exist.
In the case when the foliation by leaves is regular and the quotient is well-defined, we have a Riemannian submersion (cf.~\cite{O'Neill66}). Otherwise, as it was proven in \cite{Kotov-Strobl16_Geom1}, the obtained foliation is Riemannian (in the sense of \cite{Molino88}): every geodesic curve with an initial velocity orthogonal to the fibers of the foliation remains orthogonal to the fibers. Notice that in this context the Lie algebroid is not endowed with a fiber metric.

\vskip 3mm\noindent In \cite{Kotov-Strobl_CYMH15}, this setting was extended by precisely such a fiber metric:
the main driving motivation for this study was given by an action functional which generalizes the Yang-Mills model with Higgs fields by replacing its structure Lie group together with its action on the space of Higgs fields $M$ by a more general Lie groupoid over $M$.
The  functional of this Curved Yang-Mills-Higgs model (CYMH) is defined whenever one is given the following data: a Lorentzian manifold $(\Sigma,\gamma)$, which serves as the space-time manifold of the theory, and a Riemannian manifold $(M,g)$ together with a Lie algebroid $(A, \rho)$ over $M$, which is supplied with a fiber metric $\kappa$ and a linear connection $\nabla$, as well as
an $A$-valued $2$-form $B$ on $M$. The fields of CYMH are bundle maps from the source tangent bundle $T\Sigma$ to the target Lie algebroid $A$ viewed as a vector bundle over $M$.
The following theorem provides the compatibility conditions between the data on the target.

\begin{theorem}[Kotov-Strobl, \cite{Kotov-Strobl_CYMH15}]\label{thm:CYMH}
 The CYMH-functional $S_{CYMH}$  \cite{Kotov-Strobl_CYMH15} is gauge invariant
 if and only if the following conditions hold:
 \begin{enumerate}
  \item $\nabla$ is a Cartan connection on the target Lie algebroid $A$
  \item ${}^\tau{}\nabla (g)=0$;
  \item ${}^\a{}\nabla (\kappa)=0$;
  \item $R_{\nabla}+[\mathrm{D},\iota_\rho] (B)+\langle t, B\rangle =0$. \end{enumerate}
 \end{theorem}

\noindent Here we use the notations from \cite{Kotov-Strobl16_Geom1}
for the Lie algebroid connections ${}^\tau{}\nabla$ and ${}^\a{}\nabla$, induced by the ordinary connection $\nabla$ on $A$
(cf.~\cite{Blaom04,Blaom05} as well as \cite{BKS04,Strobl_AYM04,Mayer-Strobl09} for the related gauge transformations). In the last condition above, $R_\nabla$ is the curvature and $\mathrm{D}$  the exterior covariant derivative of $\nabla$,  $t$ is the $A$-torsion of  ${}^\alpha{}\nabla$, and
$\iota_\rho$ denotes the contraction with the $TM$-part of $\rho \in \Gamma(A^* \otimes TM)$.
The first condition in Theorem \ref{thm:CYMH} requires that the connection  $\nabla$ on $A$ is compatible with the Lie algebroid structure of $A$; it turns  $A$ into a Cartan-Lie algebroid (see Defintion \ref{def:Cartan-Lie_alg} below and \cite{Kotov-Strobl16_Geom1,Blaom04,Blaom05} for more details). The second condition expresses that $A$ is a Killing Lie algebroid
over $(M,g)$, as explained above. Likewise, we will call a Lie algebroid $A$ together with a fiber metric $\kappa$ and a vector bundle connection $\nabla$ satisfying the equation $(3)$ a metric Lie algebroid.

\vskip 3mm\noindent In contrast to the first two conditions, the third and the fourth ones were not yet studied in detail.

\vskip 3mm\noindent In the present paper, we give a geometric interpretation of the first three conditions in terms of a reflection- and bi-invariant  metric
on a Cartan-Lie groupoid, turning the latter one into a particular example of a Riemannian groupoid in the sense of  \cite{GGHR1989}; we call such Lie groupoids Riemannian Cartan-Lie groupoids or simply Riemann-Cartan groupoids and their Lie algebroids, which are canonically equipped with the data $(\kappa,g)$ satisfying the conditions (1), (2), and (3) above, quadratic Lie algebroids; in fact, in this case even positive quadratic Lie algebroids since $g$ and $\kappa$ are definite, while the notion generalising quadratic Lie algebras evidently does not need this restriction. Riemannian groupoides in the sense of \cite{delHoyo-Fernandes15} require additional smooth metrics on groupoid nerves; at the end of this paper we provide a sufficient condition for such an extension to exist.

\vskip 3mm\noindent The paper is organized as follows. In Section \ref{sec:Jet_groupoids}
 we summarize the relevant material on jet groupoids. In Section \ref{sec:Cartan-Lie_groupoids} we focus on Cartan Lie groupoids.
In the last Section \ref{sec:bi-invariant_metrics}, we prove the main theorem of this paper: We provide the necessary and sufficient conditions for that
 a Cartan-Lie groupoid $\calG$ whose Lie algebroid $A$ is equipped with a $\calG$-invariant positive quadratic structure\footnote{In the case that $\calG$ is source-connected, the Ad$_\calG$-invariance follows directly from the infinitesimal ad-invariance built into the definition of a quadratic Lie algebroid. We assume that the Cartan structure on $A$ is the canonical one induced from $\calG$ here. Cf.~also Proposition \ref{prop:Cartan-Lie_algebroid_vs_groupoid} and the remark following Definition \ref{def:quadratic_Lie_algebroids} in the main text below.}
 can be endowed with a bi-invariant and inversion invariant metric $\eta$ such that its restriction to $\ker(\md t)\vert_M\cong A$
 coincides with $\kappa$ and the source map $s$ is a Riemannian submersion onto $(M,g)$.

\vskip 3mm \noindent A.K. would like to thank Projeto P.V.E. 88881.030367/2013-01 (CAPES/Brazil),
FAPESP grant 2011/11973-4 for funding his visit to ICTP-SAIFR in December 2016 and April 2017 where part of this work was done,
and the institutional support of UFPR (Curitiba) and the University Hradec Kralove. T.S. is grateful to the CNRS for according a delegation and an attachment to beautiful IMPA for one semester and equally to IMPA and its staff for their hospitality during this stay.



\section{Jet groupoids}\label{sec:Jet_groupoids}

\setcounter{theorem}{0}

\noindent In this section we restrict our attention to
jet groupoids. A more general discussion of jet spaces can be found in \cite{KV11,Kruglikov-Lychagin08} (see also \cite{KLV86,Lychagin95,Olver86} and the references given therein).
For jet groupoids and algebroids we refer to \cite{Kumpera-Spencer1972, Kumpera1975}.
We presume the knowledge of the definition and the most basic properties of Lie groupoids (see for instance \cite{Mackenzie05} or \cite{Crainic-Fernandes2011}).

\vskip 3mm\noindent Let $\calG$ be a Lie groupoid over a base $M$ whose source and target maps are $s$ and $t$,
respectively. Let $1_x$ be the identity at $x\in M$ and $\iota$ the inversion map.
Recall that the multiplication operation on a groupoid $\calG$ is a smooth map $m\colon \calG^2\to \calG$, where $\calG^2$ is the fibered product\footnote{The fibered product $N_{f}\!\times_{f'}N'$  of $f\colon N\to M$ and $f'\colon N'\to M$
is a smooth manifold, if at least one of the maps is a surjective submersion.
We shall also use the notation $N\times_M N'$ for the fibered product over $M$ whenever all maps are clearly defined.} over $M$
\beqn \calG^2=\calG_s\!\times_t\calG=\{(u,u')\in \calG\times\calG \, \mid \, s(u)=t(u')\}\,\eeq
which meets the criteria of associativity. Hereafter we denote the product $m (u,u')$
by $u* u'$ whenever the elements $u$ and $u'$ are composable.


\begin{deff}\label{deff:G-space} Let $Z$ be a smooth manifold. A left $\calG$-action on $Z$ is a pair
$(\mu, a)$ of maps $\mu\colon Z\to M$ and $a\colon \calG_s\!\times_\mu Z\to Z$ such that $\mu (u z)=t(u)$
for all $(u,z)\in \calG_s\!\times_\mu Z$,
i.e.~the following diagram is commutative
\bqa\nonumber
\xymatrix{ \calG_s\!\times_\mu Z \ar[rr]^a \ar[drr]_{t\circ\mathrm{pr}_1}  && Z \ar[d]^{\mu \qquad ,} \\
   && M  }
\eqa
and the action satisfies
$u(u' z)=(u*u')z$ and $1_{\mu(z)}z=z$ for all compatible $u,u'\in\calG$ and $z\in Z$. Here $\mathrm{pr}_i$ is the projection of the fibered product
onto the corresponding factor and $uz$ is a short notation for $a (u,z)$.
We shall call $Z$ a (left) $\calG$-space and $\mu$ a moment map.
A right $\calG$-action and a right $\calG$-space are defined in a similar way.
\end{deff}

\begin{deff}\label{def:morphism_G-spaces} A smooth map $f$ between $\calG$-spaces
$(Z,\mu,a)$ and $(Z',\mu',a')$ is a morphism of $\calG$-spaces, if it commutes with the corresponding structure maps, i.e.~if
the following diagrams are commutative.
\bqa\nonumber
\xymatrix{ Z \ar[rr]^f \ar[dr]_{\mu}  && Z' \ar[dl]^{\mu'} \\
   & M &   }
\hspace{20mm}
\xymatrix{ \calG_s\!\times_\mu Z \ar[rr]^{\Id\times f} \ar[d]_{a} && \calG_s\!\times_\mu Z' \ar[d]^{a'} \\
Z \ar[rr]^f   && Z'  }
\eqa
$\calG$-spaces together with their morphisms form a category.
\end{deff}

\noindent Hereafter we shall assume, for simplicity, that $\mu$ is always a surjective submersion.

\begin{deff}\label{deff:G-equivariant_bundle} A smooth  bundle $\pi\colon E\to Z$ over a $\calG$-space $Z$ is called $\calG$-equivariant,
if $E$ is a $\calG$-space and $\pi$ is a morphism of $\calG$-spaces. A section of a $\calG$-equivariant bundle is called $\calG$-invariant, if it is
a morphism of the corresponding $\calG$-spaces.
\end{deff}

\begin{example}\label{ex:base_as_a_groupoid_space} The base $M$ of a groupoid $\calG$ is a left and right $\calG$-space simultaneously, where $\mu=\Id$ and $a=t\circ \mathrm{pr}_1$ or
$s\circ \mathrm{pr}_2$
depending on whether we view $M$ as left or right $\calG$-space.
\end{example}

\begin{example}\label{ex:groupoid_as_a_groupoid_space} $(\calG,t,m)$ and $(\calG,s,m)$ are canonical left and right $\calG$-spaces, respectively; furthermore, they can be thought of as
$\calG$-equivariant bundles over $M$.
\end{example}


\noindent Just like local diffeomorphisms of a smooth manifold, local bisections of a Lie groupoid
are endowed with a (partially defined) operation of composition, provided their domains are compatible.  In the present paper we define a local bisection  $\Sigma$ of a groupoid $\calG$ to be
a submanifold of $\calG$ such that $s(\Sigma)$ and $t(\Sigma)$ are open subsets of $M$ and $s\!\!\mid_{\Sigma}\colon \Sigma\to s(\Sigma)$ and $t\!\!\mid_{\Sigma}\colon \Sigma\to t(\Sigma)$ are diffeomorphisms. In what follows we abbreviate  $s\!\!\mid_{\Sigma}$ by $s_{\sst\Sigma}$ and likewise $t\!\!\mid_{\Sigma}$ by $t_{\sst\Sigma}$.
Let $\Sigma$ and $\Sigma'$ be local bisections. Their product is defined iff
$s(\Sigma) \cap t(\Sigma') \neq	\emptyset$, in which case it equals to
\beq\label{product_of_bisections}
\Sigma *\Sigma' = m\left(\Sigma\, _s\!\times_t\Sigma'\right)\, ,
\eeq
where $\Sigma\, _s\!\times_t\Sigma' = \{(q, q')\in\Sigma\times \Sigma' \, \mid  \, s(q)=t(q')\}\subset\calG^2$.
In other words, their product 
is
\beqn
\Sigma *\Sigma' =\{q* q' \, \mid \, q\in\Sigma, \, q'\in\Sigma', \, s(q)=t(q')\}\,.
\eeq

\noindent While the space $\cB (\calG)$ of global bisections of a Lie groupoid is a group, the space $\cB_{\sst loc} (\calG)$ of local bisection is a (generalized) pseudogroup (in the sense of \cite{Salazar2013}).\footnote{See also http://mathworld.wolfram.com/Pseudogroup.html and the reference contained there for the definition of a pseudogroup.}
Every local bisection $\Sigma$ gives us two locally defined diffeomorphisms of $\calG$,
the left and right translations, which we denote by $L_{\sst\Sigma}$
and $R_{\sst\Sigma}$, respectively. More precisely, the left translation operator is acting from $t^{-1}s(\Sigma)$ to $t^{-1}t(\Sigma)$ by the formula
$L_{\sst\Sigma} (u)=u'*u$, where $u'$ is the unique element of $\Sigma$, such that $s(u')=t(u)$. The right translation operator
is defined in a similar way and is acting from $s^{-1}t(\Sigma)$ to $s^{-1}s(\Sigma)$.
The left translation $L_{\sst\Sigma}$ covers a diffeomorphism $\varphi_{\sst\Sigma}=t_{\sst\Sigma}s^{-1}_{\sst\Sigma}\colon s(\Sigma)\to t(\Sigma)$
with respect to the surjective submersion $t$, while the right translation $R_{\sst\Sigma}$ covers
$\varphi^{-1}_{\sst\Sigma}=s_{\sst\Sigma}t^{-1}_{\sst\Sigma}\colon t(\Sigma)\to s(\Sigma)$ with respect to the surjective submersion $s$, i.e.~the following diagrams are commutative.
\beq\label{left-right_section_translations}
\xymatrix{ t^{-1}s(\Sigma) \ar[rr]^{L_{\sst\Sigma}} \ar[d]_{t}  && t^{-1}t(\Sigma) \ar[d]_{t}\\
  s(\Sigma)\ar[rr]^{\varphi_{\sst\Sigma}} & & t(\Sigma)  }
\hspace{20mm}
\xymatrix{ s^{-1}t(\Sigma) \ar[rr]^{R_{\sst\Sigma}} \ar[d]_{s}  && s^{-1}s(\Sigma) \ar[d]_{s}\\
  t(\Sigma)\ar[rr]^{\varphi^{-1}_{\sst\Sigma}} & & s(\Sigma)  }
\eeq
It is obvious that, for any pair of bisections with compatible domains, $L_{\sst \Sigma *\Sigma'}=L_{\sst \Sigma}L_{\sst \Sigma'}$,
$\varphi_{\sst \Sigma *\Sigma'}=\varphi_{\sst \Sigma}\varphi_{\sst \Sigma'}$, and $R_{\sst \Sigma *\Sigma'}=R_{\sst \Sigma'}R_{\sst \Sigma}$. The inverse map is given by the "inverse bisection" $\iota (\Sigma)$, which satisfies in particular
$\varphi_{\sst \iota(\Sigma)}=\varphi^{-1}_{\sst \Sigma}$. Therefore the correspondence $\Sigma\mto \varphi_{\sst \Sigma}$ is a morphism
from $\cB_{\sst loc} (\calG)$ to the pseudogroup of local diffeomorphisms of $M$.

\vskip 2mm\noindent The equivariance property of the left and right translations by $\Sigma$,
expressed via the commutative diagrams (\ref{left-right_section_translations}), admits the following generalization.
 \begin{proposition}\label{prop:bisection_action}$\!\!$\footnote{We remark that we deal with
 right $\calG$-equivariant bundles over the base of a Lie groupoid in Proposition \ref{prop:bisection_action}  and equally so in Proposition \ref{prop:groupoid_action_prolongation} below. A right action is more convenient for the corresponding infinitesimal
 (Lie algebroid) action; obviously there exist similar left $\calG$-equivariant versions of the above mentioned propositions.}{\mbox{}\vskip 2mm}
 \begin{enumerate}
 \item Let $(Z,\mu,a)$ be a right $\calG$-space. Then there is a morphism from $\cB_{\sst loc} (\calG)$ to the pseudogroup of local diffeomorphisms of $Z$,
 such that each local bisection $\Sigma$ determines a diffeomorphism from $\mu^{-1}t(\Sigma)$ to $\mu^{-1}s(\Sigma)$.
 \item Let $\pi$ be a right $\calG$-equivariant bundle over $Z$. Then the pseudogroup $\cB_{\sst loc} (\calG)$ is acting by local bundle isomorphisms.
 \item Given a local section $\sigma$ of $\pi$ over $\mu^{-1}t(\Sigma)$, we obtain a new local section over $\mu^{-1}s(\Sigma)$, denoted by $\sigma\Sigma$.
 \end{enumerate}
 \end{proposition}
 \noindent\proof For each $z\in \mu^{-1}s(\Sigma)$ we define $R_{\sst\Sigma}z=zu$, where $u$ is the unique element of $\Sigma$, such that
 $s(u)=\mu(z)$. One can verify that $R_{\sst\Sigma}$ is a diffeomorphism, the inverse of which is given by $\iota(\Sigma)$, and the correspondence
 $\Sigma\mto R_{\sst\Sigma}$ obeys the requirements of a right representation. The second statement follows from the equivariance of $\pi$ and
 the third one from the general property of local bundle isomorphisms. More precisely, if $\sigma$ is a local section of $\pi$ whose domain $U$ is contained in the domain of $R_{\sst\Sigma}^{-1}$, i.e.~$\mu (U)\subset s(\Sigma)$, then $R_{\sst\Sigma}^*\sigma$ is a section of $\pi$ over $R_{\sst\Sigma}^{-1}U$, defined by the formula
 $  R_{\sst\Sigma}^*\sigma (x)= R_{\sst\Sigma}^{-1}\sigma \left(R_{\sst\Sigma}x\right)$ for all $x\in R_{\sst\Sigma}^{-1}U$.
  $\blacksquare$


\vskip 2mm\noindent Denote by $V(s)$ and $V(t)$ the subbundles of vectors that are tangent to the fibers of the corresponding projections to $M$, i.e.~$V(s)=\Ker\, \md s$,
 $V(t)=\Ker\, \md t$; henceforth we shall name them $s$-vertical and $t$-vertical, respectively.

\begin{rem}\label{ref:groupoid_does_not_act_on_tangent}
In contrast to Lie groups, left and right translations are well-defined only on $t$-vertical and $s$-vertical vectors, respectively. Therefore
$V(t)$ and $V(s)$ are left  and right $\calG$-equivariant bundles over $\calG$, respectively (see Example \ref{ex:groupoid_as_a_groupoid_space}).
\end{rem}

\noindent Let us recall that sections of the Lie algebroid $A$ of $\calG$ are in one-to-one correspondence with $t$-vertical left $\calG$-invariant vector fields on $\calG$.
  The Lie bracket on sections of $A$ is induced by the Lie bracket of the corresponding vector fields on $\calG$, while the
   anchor map $\rho$ is determined by $\md s$.


\begin{rem}\label{rem:projectible_vector_fields}
By the construction of the Lie algebroid $A$ of $\calG$,
there is a canonical $C^\infty (M)$-linear embedding $\Gamma (A)\hookrightarrow \Gamma (T\calG)^s$, i.e.~into the $s$-projectable vector fields on $\calG$,
 which respects the Lie brackets. In other words, $A$ is canonically acting on $s\colon\calG\to M$ (cf.~\cite{Higgins-Mackenzie90,Kosmann-Mackenzie02} for the definition of a Lie algebroid action).
 \end{rem}

\noindent Let $\xi$ be a section of the Lie algebroid $A$ regarded as a left-invariant $t$-vertical vector field on $\calG$,
which we shall denote by the same letter, and $\psi_\e$ the flow on $\calG$ generated by $\xi$.

\begin{rem}\label{rem:flow_in_non-compact_case} In the non-compact case
we can not guarantee that $\psi_\e$ is globally defined; it should be thought of as a family of local diffeomorphisms $\{\psi^\a_\e\}$,
solving the corresponding Cauchy problem over their domains $U^\a$ with the initial condition 
$\psi^\a_\e (u)|_{\e=0}=u$ for each $u\in U^\a$,
such that $\{U^\a\}$ gives us an open cover of $\calG$. Although the time intervals depend on the index $\a$, one has
$\psi^\a_\e =\psi^\beta_\e$ over $U^\a\cap U^\b$ at a fixed time $\e$ from the intersection of the corresponding time intervals.
\end{rem}

\noindent From $\md t(\xi)=0$ we conclude that $t\circ \psi_\e =t$. Taking into account that $\psi_\e$ is left $\calG$-invariant, we deduce that
the flow of $\xi$ is given by the right multiplication on a $1$-parameter family of $t$-sections $\la_\e$, which are also bisections,
where $\lambda_\e (x)=\psi_\e (1_x)$ for $x\in M$.
This family $\la_\e$ covers the flow\footnote{In the non-compact case we have a family of local bisections
since the flow of $\xi$ on $\calG$ as well as the flow of $\rho(\xi)$ on $M$ may not be globally defined,
see Remark \ref{rem:flow_in_non-compact_case}. However, such a local description is sufficient to determine the infinitesimal
action of the Lie algebroid of $\calG$ whenever we have a representation of the pseudogroup $\cB_{\sst loc} (\calG)$ by local diffeomorphisms;
 we employ this point of view in Proposition \ref{prop:LA_action_vs_equiv_bundle}.} of $\rho(\xi)$ on $M$ (see the right diagram in (\ref{left-right_section_translations})).
As a corollary of Proposition \ref{prop:bisection_action},
we get the following relation between $\calG$-equivariant bundles
(Definition \ref{deff:G-equivariant_bundle}) and Lie algebroid actions.

\begin{proposition}\label{prop:LA_action_vs_equiv_bundle} Let $\pi\colon E\to M$ be a right $\calG$-equivariant bundle, then
the corresponding infinitesimal flow gives us a structure of an infinitesimal Lie algebroid action.
\end{proposition}


\vskip 3mm\noindent
Let $\Sigma\subset\calG$ be a local bisection  such that $u\in\Sigma$. We denote by $[\Sigma]^{k}_u$ the $k$-jet of $\Sigma$ at $u$.

\begin{deff}\label{def:k-groupoid}
 Given a Lie groupoid $\calG$, for every $k\ge 0$, we denote by $J^k (\calG)$ the space of $k$-th order jets of local bisections of $\calG$.
 \end{deff}

\begin{rem}\label{rem:bisection_jets_inside_jets}
Taking into account that every (local) bisection determines a (local) section of both $s\colon\calG\to M$ and $t\colon\calG\to M$, we can identify
$J^k(\calG)$ with an (open dense) subset of $J^k (s)$ and $J^k (t)$, simultaneously, where, by definition, $J^k (\pi)$ is the space of $k-$jets of local sections of a bundle $\pi$.
Moreover, we have $J^k(\calG)=J^k (s)\cap J^k (t)$, provided
the jet spaces of sections are viewed as (open dense) subsets of the space of jets of $n$-dimensional submanifolds of $\calG$, where $n=\dim M$.
\end{rem}

\noindent
The multiplication of local bisections with compatible domains induces a multiplication law on $J^k (\calG)$ which turns  $J^k (\calG)$  into a Lie groupoid,
such that the natural projection maps $\pi_{\sst k,l}\colon J^k (\calG)\to J^l(\calG)$ for all $l\le k$ are Lie groupoid morphisms.
More precisely, let $u\in\Sigma$, $u'\in\Sigma'$ such that $s(u)=t(u')$, then $[\Sigma]^{k}_u*[\Sigma']^{k}_{u'}\colon=[\Sigma *\Sigma']^{ k}_{u * u'}$.
 One notes that the multiplication is smooth and the so-defined product depends only on the $k$-jets  $[\Sigma]^{ k}_u$ and $[\Sigma']^{ k}_{u'}$
 and not on the choice of representatives
$\Sigma$ and $\Sigma'$.
 In particular, if we denote the source and target maps of $J^k (\calG)$ by $s_{\sst k}$ and $t_{\sst k}$, respectively, such that $s=s_{\sst 0}$ and $t=t_{\sst 0}$,
 we obtain the equalities $s_{\sst k}=s_{\sst l}\circ \pi_{\sst k,l}$
 and  $t_{\sst k}=t_{\sst l}\circ \pi_{\sst k,l}$. The identity at $x\in M$ is given by the $k$-jet of the identity section
 and the inverse of $[\Sigma]^k_u$ is given by $[\iota(\Sigma)]^k_{\iota(u)}$.

 \vskip 2mm\noindent Besides being Lie groupoid morphisms, the canonical projection maps $\pi_{\sst k,l}\colon J^k(\calG)\to
J^l (\calG)$ constitute an inverse system, i.e.~$\pi_{\sst k,l}=\pi_{\sst m,l}\circ \pi_{\sst k,m}$ whenever $l\le m\le k$. Thus the inverse limit
 $J^\infty (\calG)=\varprojlim J^k (\calG)$, the space of infinite jets of local bisections, is a groupoid again.

\vskip 2mm\noindent
The next proposition is about jet prolongations of a $\calG$-equivariant bundle over the base of a Lie groupoid.
 \begin{proposition}\label{prop:groupoid_action_prolongation}
 Let $\pi$ be a right $\calG$-equivariant bundle over $M$.
 Then for each $k$ from $1$ to $\infty$ there exists a canonical structure of a right $J^k(\calG)$-equivariant
 bundle on $\pi_{\sst k}$, which is uniquely determined as follows. Let $u\in\calG$, $x\in M$
 such that $s(u)=x$, let $\sigma$ be a local section of $\pi$ at $x$ and $\Sigma$ a local bisection of $\calG$ through $u$,
 such that $\sigma$ and $\iota(\Sigma)$ are compatible in the sense of Proposition \ref{prop:bisection_action}. Then one has
  $ [\sigma]_x^k[\Sigma]_u^k=[R_{\sst\iota(\Sigma)}^*\sigma]_{s(u)}^k$.
\end{proposition}
\noindent\proof The proof is straightforward; it is sufficient to show that the $k$-jet of
$ R_{\sst\iota(\Sigma)}^*\sigma$ at $s(u)$ depends only on the $k$-jets of $\sigma$ at $x$ and $\Sigma$ at $u$.
The groupoid action properties follow from Proposition \ref{prop:bisection_action}.
$\blacksquare$


\begin{deff}\label{def:Jet_algebroid}Let $(A,\rho,[\cdot,\cdot])$ be a Lie algebroid over $M$. The $k-$jet Lie algebroid $J^k (A)$
is the vector bundle of $k-$jets of sections of $A$, endowed with the canonical Lie bracket and the anchor map as follows:
the bracket in $J^k (A)$ is defined such that taking the Lie brackets commutes with the $k$-jet prolongation of sections,
 \bqa\label{jets-bracket} [j_k(\xi), j_k(\xi')] = j_k ([\xi,\xi'])
 \eqa
 for all sections $\xi,\xi' \in \Gamma(A)$, while its anchor $\rho_{\sst k}$ is fixed by the morphism property to obey
 \bqa \label{jet-anchor}
 \rho_{\sst k}(j_k(\xi))=\rho(\xi)\,.
 \eqa
\end{deff}

\vskip 1mm\noindent By the above construction of the Lie algebroid structure on $J^k (A)$, the projections $\pi_{\sst k,l}\colon J^k(A)\to J^l (A)$ for
all $l<k$ are Lie algebroid morphisms, thus the projective limit $J^\infty (A)=\varprojlim J^k (A)$ admits a canonical
Lie algebroid structure.

\vskip 3mm\noindent The next statement
(cf.~\cite{Kumpera-Spencer1972, Kumpera1975})  follows from propositions \ref{prop:groupoid_action_prolongation}
and \ref{prop:LA_action_vs_equiv_bundle}.

\begin{proposition}\label{prop:k-jet_algebroid} Let $A$ be the Lie algebroid of a Lie groupoid $\calG$, then
the Lie algebroid of $J^k (\calG)$ is $J^k (A)$ for all $k\in\N$ and the Lie algebroid of $J^\infty (\calG)$ is $J^\infty (A)$.
\end{proposition}

\begin{cor}\label{cor:LA_prolongation_vs_groupoid_prolongation}
The infinitesimal flow of the jet prolongation of a Lie groupoid action
gives us the jet prolongation of the corresponding Lie algebroid action.
\end{cor}

\begin{example}\label{ex:local_diffeo}
 Let $\calG$ be the pair groupoid associated to a smooth manifold $M$, $\calG=M\times M$, then $J^k (\calG)$ is the groupoid of $k$-jets of local diffeomorphisms of $M$.
 \end{example}

 \begin{example}\label{ex:0-jet_groupoid}
 Let $k=0$, then $J^0(\calG)=\calG$.
 \end{example}


\section{Cartan-Lie groupoids}\label{sec:Cartan-Lie_groupoids}

\setcounter{theorem}{0}

\noindent
The notion of Cartan-Lie algebroids   \cite{Blaom04,Blaom05} and their groupoids \cite{Blaom16} was pioneered by A.~Blaom; additional details on this subject can be found in these papers.
Here we follow our presentation of Cartan-Lie groupoids as given in \cite{Alexei_talk}.
For general multiplicative distributions on Lie groupoids, which play an important role in this context, we refer to \cite{CSS12} (see also \cite{Tang2004, Salazar2013}).

\vskip 2mm\noindent We start with properties of the $1$-jet groupoid $J^1 (\calG)$. Thereafter we will focus on bi-connections
on a Lie groupoid and on Cartan-Lie groupoids.

\vskip 2mm\noindent An element $[\Sigma]^1_u$ of $J^1 (\calG)$ can be uniquely identified with the tangent space
 $T_u\Sigma$.
 Now if we recall that $s\colon\Sigma\to M$ and $t\colon\Sigma\to M$ are local diffeomorphisms, we get linear isomorphisms $\md s\colon T_u\Sigma\iso T_{s(u)}M$
 and $\md t\colon T_u\Sigma\iso T_{t(u)}M$.
 With these identifications, the fiber $\pi^{-1}_{\sst 1,0}(u)$
  consists of all $n$-dimensional linear subspaces $L\subset T_u\calG$, $n=\dim M$, such that the restriction of $\md s$ and $\md t$ to $L$ are linear isomorphisms:
  \beqn
  \md s_u\colon L \iso T_{s(u)}M \, , \hspace{3mm} \md t_u\colon L \iso T_{t(u)}M \,.
  \eeq
  From now on, we shall call these subspaces $s$-projectible and $t$-projectible, respectively.
  Equivalently, $\pi^{-1}_{\sst 1,0}(u)$ coincides with the space of all $n$-dimensional linear subspaces of $T_u\calG$ the intersection of which with both
  $V(t)$ and $V(s)$ is zero. Thus a section of $\pi_{\sst 1,0}$ can be regarded as a smooth rank $n$ distribution $H$ on $\calG$, consisting of subspaces
  that are both $s$-projectible and $t$-projectible.

  \vskip 2mm\noindent Let $(u,u')\in\calG^2$ and let $L\in \pi^{-1}_{\sst 1,0}(u)$ and $L'\in \pi^{-1}_{\sst 1,0}(u')$. Application
  of formula (\ref{product_of_bisections}) yields that the product of $L$ and $L'$
  is equal to $\md m \left( L\, _{\md s}\!\times_{\md t} L'\right)$, where
  \beqn
  L\, _{\md s}\!\times_{\md t} L' =\{ (v,v')\in L\times L'\subset T_u\calG\times T_{u'}\calG \, \mid \, \md s(v)=\md t(v')\}\subset T_{(u,u')}\calG^2\,.
  \eeq

\vskip 2mm\noindent The proof of the following statement is similar to the proof of Proposition \ref{prop:k-jet_algebroid}.

\begin{lemma}\label{lem:1-jet_groupoid_tangent_action} Let $Z$ be a $\calG$-space, then its tangent bundle $TZ$ is a $J^1(\calG)$-space.
If $\pi$ is a $\calG$-equivariant bundle over $Z$, then $\md\pi$ is a $J^1(\calG)$-equivariant
bundle over $TZ$. In particular, $T\calG$ admits a canonical structure of both left and right $J^1 (\calG)$-spaces.
\end{lemma}



\begin{deff}\label{def:biconnection}
A section of $\pi_{1,0}$, regarded as a bundle without additional structure, will be called a bi-connection.
\end{deff}

\begin{rem}\label{rem:tangent_action} In a general, a Lie groupoid does not act on its tangent space, see Remark \ref{ref:groupoid_does_not_act_on_tangent}. However,
as soon as we fix a section $\sigma$ of $\pi_{1,0}$, left and right translations become well-defined on all of $T\calG$
by combining  the map $\sigma$ with the action of $J^1 (\calG)$ on $T\calG$ (see Lemma \ref{lem:1-jet_groupoid_tangent_action}),
although these translations may not be necessarily groupoid actions, unless $\sigma$ respects also the groupoid structure of the bundle,
cf.~Proposition \ref{prop:morphism_vs_multiplicativity} below.
\end{rem}

  \noindent In more details,
  let $H$ be a distribution on $\calG$ which represents a section of $\pi_{\sst 1,0}$. Let $(u,u')\in\calG^2$, $v\in T_{u}\calG$, and $v'\in T_{u'}\calG$.
  Then the left translation $L_u (v')$ is defined as $\md m (w, v')$, where $w$ is the unique vector in $H_u$ that satisfies
  $\md s(w)=\md t(v')$; likewise, the right translation $R_{u'} (v)$ is defined as $\md m (v, w')$,
  where $w'$ is the unique vector in $H_{u'}$ that satisfies
  $\md s(v)=\md t(w')$. It is obvious that the left $\calG$-action on $V(t)$ and the right $\calG$-action on $V(s)$ do not depend on the choice of a bi-connection
  and thus coincide with the canonical groupoid actions.

\begin{deff}\label{def:multiplicative}
A multiplicative distribution on $\calG$ is a smooth vector subbundle $H$ of $T\calG$, that obeys the multiplicative rule with respect to
$m$, i.e.~that satisfies $\md m \left(H\,_{\md s}\!\times_{\md t} H\right)\subset H$, where
\beqn
\left( H\,_{\md s}\!\times_{\md t} H\right)_{(u,u')}=\{(v, v')\in H_u\times H_{u'}\, \mid \, \md s (v)=\md t (v')\}\subset T_{(u,u')}\calG^2\,, \hspace{3mm} (u,u')\in \calG^2\,.
\eeq
\end{deff}

\begin{proposition}\label{prop:morphism_vs_multiplicativity}
The following properties are equivalent:
\begin{enumerate}
 \item a section of $\pi_{1,0}$ is a groupoid morphism;
 \item the corresponding distribution $H$ is multiplicative;
 \item the left and right translations determined by $H$ are left and right groupoid actions of $\calG$,
 respectively.
\end{enumerate}
\end{proposition}
\noindent\proof
The proof of this proposition follows from the explicit formula of the product in $J^1 (\calG)$
and Lemma \ref{lem:1-jet_groupoid_tangent_action}. In particular, one can verify that
the multiplicative property leads to $H|_M=TM$, where $M$ is identified, as usual, with the identity section,
and $H_{\iota (u)}=\md\iota \left(H\right)_{\iota (u)}$.
We leave the details to the reader.
$\blacksquare$

  \begin{deff}\label{def:Cartan-Lie_gpd}
 By a Cartan-Lie groupoid (or simply a Cartan groupoid) we mean a Lie groupoid together with a bi-connection that induces a groupoid morphism from ${\calG}$ to $J^1 (\calG)$.
 Equivalently, by Proposition \ref{prop:morphism_vs_multiplicativity},
 a Cartan-Lie groupoid over an $n$-dimensitonal base is a Lie groupoid supplied with
 a multiplicative rank $n$ distribution, which consists of subspaces
  that are $s$-projectible and $t$-projectible, simultaneously.
  \end{deff}

\begin{lemma}\label{lem:bi_invariance} Let $\calG$ be a Cartan groupoid and $H$ its multiplicative distribution. Then
\begin{enumerate}
\item $H$, $V(s)$, and $V(t)$ are bi-invariant under the corresponding $\calG$ action;
\item The horizontal $s$- and $t$-lifts of a base vector field are left- and right-invariant vector fields on the groupoid, respectively.
\end{enumerate}
\end{lemma}

\noindent\proof
The bi-invariance of $H$ follows from the explicit formula of the tangent groupoid action, cf.~Remark \ref{rem:tangent_action}, and the multiplicativity of $H$, Definition \ref{def:multiplicative}. 

\vskip 3mm\noindent
Now we shall prove the other statements.
Let $X$ be a vector field on $M$, $\tilde X$ be an arbitrary $s$-lift of $X$, i.e.~a vector field on $\calG$
such that $\md s(\tilde X)=X$. By the formula of left translation
$ L_u \left(\tilde X_{u'}\right)=\md m \left( w,\tilde X_{u'}\right)$,
where $(u,u')\in\calG^{2}$ and $w$ is the unique horizontal vector at $u$ such that $\md s(w)=\md t (\tilde X_{u'})$.
Denote by $\mathrm{pr}_1$ and $\mathrm{pr}_2$ the canonical projection of $\calG^{2}$ to the first and second factor, respectively. One
has $s\circ m=s\circ \mathrm{pr}_2$, thus
\beq\label{ds_of_left_action}
\md s \left( L_u \left(\tilde X_{u'}\right)\right)=\md \left( s\circ m\right) \left(w, \tilde X_{u'}\right)=\md s (\tilde X_{u'})= X_{s(u')}\,.
\eeq
In particular, if $X=0$ and thus $\tilde X\in V(s)$, then $L_u \left(\tilde X_{u'}\right)\in V_{u*u'} (s)$ for any pair $(u,u')\in\calG^{2}$. This shows that
$V(s)$ is a left $\calG$-invariant subbundle. On the other hand, $V(s)$ is also a right $\calG$-invariant subbundle, such that the right action of $\calG$ on $V(s)$
is independent of the choice of $H$, see Remark \ref{ref:groupoid_does_not_act_on_tangent}, therefore $V(s)$ is a bi-invariant subbundle.
Similarly, also $V(t)$ is a bi-invariant subbundle of $T\calG$.

\vskip 2mm\noindent
Let us prove now that the horizontal $s$-lift of a base vector field is left $\calG$-invariant.
Let $\tilde X^s$ be the $H$-horizontal $s$-lift of $X$, i.e.~the unique section of $H$
such that $\md s(\tilde X^s)=X$. Let $(u,u')\in\calG^2$, then using (\ref{ds_of_left_action})
and the multiplicativity of $H$ and the uniqueness of the horizontal $s$-lift of a base vector field, we obtain $L_u \left(\tilde X^s_{u'}\right)=\tilde X^s_{u*u'}$. By a similar argument,
we show that the horizontal $t$-lift of a base vector field is $\calG$-right invariant. This completes the proof.
$\blacksquare$


 \begin{deff}\label{def:Cartan-Lie_alg}  $(A,\nabla)$  is called a Cartan Lie algebroid over $M$, if $A$ is a Lie algebroid, $\nabla$ a connection on $A\to M$, and the induced splitting of the Bott sequence $\sigma \colon A \to J^1(A)$ \cite{Bott} is a Lie algebroid morphism. The corresponding connection $\nabla$ is called a Cartan connection on $A$.
 \end{deff}

\begin{proposition}\label{prop:Cartan-Lie_algebroid_vs_groupoid}
The Lie algebroid of a Cartan groupoid admits a canonical Cartan structure. Vise-versa,
a Cartan structure on the Lie algebroid of a target simply-connected Lie groupoid can be integrated to a Cartan structure on this groupoid.\end{proposition}

\noindent\proof By definition, a Cartan groupoid is a Lie groupoid $\calG$ together
with a section of $\pi_{\sst 1,0}\colon J^1 (\calG)\to \calG$, which is a Lie groupoid morphism.
The Lie algebroid of $J^1 (\calG)$ is $J^1 (A)$, where $A$ is the Lie algebroid of $\calG$
(Proposition \ref{prop:k-jet_algebroid} for $k=1$). Thus we obtain a canonical
section of  $\pi_{\sst 1,0}\colon J^1 (A)\to A$, which is a Lie algebroid morphism. Consequently, $A$ is
a Cartan-Lie algebroid according to Definition \ref{def:Cartan-Lie_alg}.

\vskip 2mm\noindent On the other hand, a Cartan structure on the Lie algebroid $A$ gives us an action of $A$ on $T\calG$
and $J^1 (\calG)$, commuting with the left $J^1 (\calG)$ actions on $T\calG$
and $J^1 (\calG)$, respectively. This action determines a left $J^1 (\calG)$-invariant vector bundle connection on $T\calG$
along $t$-fibers. If all $t$-fibers are simply-connected, then the corresponding parallel transport along $t$-fibers, which is also left $J^1 (\calG)$-invariant,
gives us an action of $\calG$ on both $T\calG$ and $J^1 (\calG)$, commuting with the left $J^1 (\calG)$ action.
We obtain the required multiplicative distribution by applying the above mentioned parallel transport to $TM$, regarded as
the tangent bundle of the identity section of $\calG$.
$\blacksquare$

\begin{deff}\label{def:invariant_tensors}
 A tensor $\Psi$ on a Cartan groupoid $\calG$ is called left- (right-) invariant, if it is invariant under left (right) translations;
 $\Psi$ is called bi-invariant, if it is invariant under both left and right translations.
\end{deff}

\begin{deff}\label{def:adjoint_action}
The adjoint action of a Cartan groupoid $\calG$ on $\left(T\calG|_M\right)$ is
defined as follows:
\beq\label{adjoint_action}
Ad_u (v) =L_u^{-1} R_u (v)\, , \hspace{3mm} v\in T_x\calG\, ,  \hspace{1mm} u\in\calG\, , \hspace{1mm} s(u)=x\,.
\eeq
\end{deff}

\begin{lemma}\label{lem:invariant_tensors}
Every section $\Psi_0$ of $\left(T\calG|_M\right)^{\otimes p}\otimes \left(T^*\calG|_M\right)^{\otimes q}$, where $M$
is identified with the identity section, can be uniquely extended by left- (right-) translations to a left- (right-) invariant
$(p,q)$-tensor on $\calG$. This extension is bi-invariant, if and only if $\Psi_0$ is invariant under the corresponding tensor power
of the adjoint action (\ref{adjoint_action}).
\end{lemma}
\noindent\proof Straightforward generalization of the similar statement known for Lie groups.
$\blacksquare$

\begin{lemma}\label{lem:adjoint_action} The restriction of $T\calG$ to the identity section
splits into the direct sum of $Ad_\calG$-invariant subbundles in two ways:
\beq\label{decomposition_of_TG} T\calG|_M = V(t)|_M\oplus TM = V(s)|_M\oplus TM\,.
\eeq
The corresponding Lie algebroid
action on $A\simeq V(t)|_M$ and $TM$ coincides with ${}^\a{}\nabla$ and ${}^\tau{}\nabla$, respectively.
Here  ${}^\a{}\nabla$ and ${}^\tau{}\nabla$ are the representations of $A$ on $A$ and $TM$, obtained by combining of
the corresponding canonical representations of $J^1$ with a splitting of the Bott exact sequence for
$J^1 (A)$, determined by the Cartan connection.\footnote{We use the notation from \cite{Kotov-Strobl16_Geom1}.
The original theory belongs to A. Blaom (\cite{Blaom04,Blaom05}). The special case of the representation ${}^\alpha{}\nabla$ was found in \cite{Strobl_AYM04,Mayer-Strobl09} independently.}
\end{lemma}
\noindent\proof The $Ad_\calG$-invariance of $V(t)|_M$, $V(s)|_M$, and $TM$ follows from Lemma \ref{lem:bi_invariance}. By the functorial property of the correspondence between
Cartan-Lie algebroids and Cartan-Lie groupoids, the induced Lie algebroid action
on $A$ and $TM$ is necessarily ${}^\a{}\nabla$ and ${}^\tau{}\nabla$, respectively.
$\blacksquare$


\section{Riemannian Cartan-Lie groupoids}\label{sec:bi-invariant_metrics}

\setcounter{theorem}{0}

\begin{lemma}
Let $\calG$ be a Cartan-Lie groupoid and let $\eta$ be
a bi-invariant Riemannian metric on $\calG$, that is, for every pair of composable elements $u_1$ and $u_2$ of $\calG$ one has
$L_{u_1}^* \left(\eta_{u_1*u_2}\right)=\eta_{u_2}$
and $R_{u_2}^*\left(\eta_{u_1*u_2}\right)=\eta_{u_1}$. Then there exist metrics $g_1$ and $g_2$ on $M$ such that
$s\colon (\calG,\eta)\to (M,g_1)$ and $t\colon (\calG,\eta)\to (M,g_2)$ are Riemannian submersions.
If, in addition, $\eta$ is preserved by the inversion map, then necessarily $g_1=g_2$.
\end{lemma}
\noindent\proof
The orbits of the left- (right-) $\calG$-action on itself are precisely the source (target) fibers and this induces a canonical $\calG$-action on the respective conormal bundles, called the transversal action, which then is defined already without the choice of a biconnection. Now,
the left- (right-) $\calG$-invariance of $\eta$ implies that $\eta$ is transversally invariant under these left- (right-) translations. Therefore
$s$ and $t$ are Riemannian submersions for an appropriative choice of Riemannian metrics on the base. Given that the inversion map $\iota$
interchanges source and target fibers, $\iota$-invariance implies $g_1=g_2$.
$\blacksquare$

\begin{deff}\label{deff:Cartan-Lie_Riemannian_groupoid}
A Cartan groupoid together with a bi-invariant metric which is invariant under the inversion
map will be called a Riemannian Carton-Lie groupoid or also a Riemann-Cartan groupoid.
\end{deff}

\begin{rem}
If $\eta$ is invariant under the inversion and $s$ is a Riemannian submersion, then
$t$ is automatically a Riemannian submersion onto the same Riemannian base.
\end{rem}

\begin{deff}\label{def:quadratic_Lie_algebroids}
Let $A$ be a Cartan-Lie algebroid over $M$, $g$ and $\kappa$ be non-degenerate metrics on $TM$ and $A$ invariant under the representations  ${}^\tau{}\nabla$ and ${}^\a{}\nabla$, respectively. Then $(A, g, \kappa)$ will be
called a quadratic Lie algebroid. If $g$ and $\kappa$ are both positive-definite then $(A, g, \kappa)$ will be
called a positive quadratic Lie algebroid or, equivalently, we will say that $A$ is endowed with a positive quadratic structure.
\end{deff}

\vskip 1mm\noindent Provided a Cartan groupoid is at least source-connected, every metric on $TM$ or $A$, which is invariant under the representations  ${}^\tau{}\nabla$ or ${}^\a{}\nabla$,
respectively, will be also invariant under the corresponding adjoint actions of $\calG$.

\begin{deff}\label{def:compatible_CLR_metric} Let $\calG$ be a Cartan groupoid over $M$, the Lie algebroid of which is $A$.
Let $(g,\kappa)$ be a positive quadratic structure on $A$ and $\eta$ be
a bi-invariant metric on $T\calG$ turning $\calG$ into a Riemann-Cartan groupoid. We shall say that $\eta$ is compatible with $(g,\kappa)$
if $s\colon (\calG,\eta)\to (M,g)$ is a Riemannian submersion and the restriction of $\eta$ to the fibers to $A$ coincides with $\kappa$;
here $A$ is identified with $V(t)|_M$, the restriction of $\Ker \,\md t$ onto the identity section of $\calG$.
\end{deff}

\vskip 1mm\noindent The next theorem, inspired by Theorem \ref{thm:CYMH},
 is about the existence and uniqueness of a Riemann-Cartan structure on a Cartan groupoid $\calG$ compatible with a
given positive quadratic Lie algebroid structure on $A$.

\begin{theorem}\label{thm:bi-invariant_metric}{\mbox{}\vskip 2mm}
\begin{enumerate}
\item
Let $(\calG,\eta)$ be a Riemann-Cartan groupoid over $M$, $A$ its Lie algebroid, and
$(g,\kappa)$  a positive quadratic structure on $A$ compatible with $\eta$. Denote by $\rho^*$ the conjugate  of $\rho\colon A \to TM$ with respect to $g$ and $\kappa$. Then
\begin{equation}�\label{smooth}
\sqrt{1-\rho\rho^*}\in \Gamma(\End(TM))\, ,
\end{equation}
i.e.~one has the bound $\rho \rho*\le 1$ and the square root above is smooth.
\item Let $\calG$ be a Cartan groupoid over a connected base $M$, $A$ its Cartan-Lie algebroid, and $(g,\kappa)$ a positive quadratic structure on $A$,
 invariant under the adjoint action of $\calG$, such that \eqref{smooth} holds true. If either $A$ is a non-transitive Lie algebroid or $\rho_x\rho_x^*=1$ for at least one $x\in M$,
then there exists a unique Riemann-Cartan structure  $\eta$ on $\calG$ compatible with $(g,\kappa)$. Otherwise there exist precisely two such compatible metrics.
 \end{enumerate}
\end{theorem}


\noindent\proof By Lemma \ref{lem:invariant_tensors}, every bi-invariant metric $\eta$ is uniquely determined by its restriction $\eta_0$ to
$T\calG|_M$; this restriction has to be $Ad_\calG$-invariant. Hereafter, for our convenience,
we shall use the simplified notations $V(s)$ and $V(t)$ for the restrictions of the corresponding subbundles of $T\calG$ onto the identity section.

\vskip 2mm\noindent Due to the bi-invariance, $\eta$ gives us a Riemann-Cartan structure on $\calG$ compatible with $(g,\kappa)$ if and only if
 \begin{enumerate}[label=(\alph*)]
 \item the projection $\md s\colon (T\calG|_M,\eta_0)\to (TM,g)$
induces, for all $x\in M$, an isometry $\md s\colon V_x(s)^\perp\to T_x M$ or, equivalently,
$\md s^*\colon T_xM \to T\calG|_x$  where $*$  is the conjugation with respect to the corresponding
metrics;
\item the restriction of $\eta_0$ onto $V(t)$ coincides with $\kappa$ (as usual, we identify $V(t)$ with the Lie algebroid $A$);
 \item $I=\md\iota|_M$ is an isometry of $(T\calG|_M,\eta_0)$.
\end{enumerate}

\vskip 2mm\noindent Combining conditions (b) and (c) and taking into account that $V(s)=I\left(V(t)\right)$, we conclude that the metric on $V(s)$ is uniquely determined by $\kappa$.
Now we have an exact sequence of metric vector bundles $$0\to V(s)\to T\calG|_M \xrightarrow{\md s} TM\to 0\,,$$ where the second arrow as well as $\md s^*$ are fiberwise linear isometries.
The metrics on the image and kernel of $\md s$ do not uniquely determine the metric $\eta_0$ on $T\calG|_M$, unless we know the orthogonal complement to $V(s)$ in $T\calG|_M$ with respect to $\eta_0$.
 Therefore metrics on $T\calG|_M$, which fit into the above exact sequence, are in one-to-one correspondence
with complements to $V(s)$. The latter can be uniquely represented as $I V(t)^\perp$, where
$ V(t)^\perp$ is some complement to $V(t)$.
In turn, $ V(t)^\perp$ can be uniquely written as the graph of a
bundle map $$\Psi\colon TM\to A\, .$$ From now the problem is to find a "suitable" $\Psi$, such that the corresponding metric on $T\calG|_M$
satisfies the conditions (a-c).

\vskip 2mm\noindent Let $v\in T_xM$ and put  $v=v^\perp +v^s$, where $v^\perp =v+I\Psi v\in V_x(s)^\perp$ and $v^s=-I\Psi v\in V_x(s)$. Since $\md s(v)=v$, we obtain
\beqn
\eta_0 (v,v)=\eta_0 \left( v^\perp, v^\perp\right)+\eta_0 \left(v^s, v^s\right)=g(v,v)+\kappa(\Psi v, \Psi v)=g \left(v, (1+\Psi^*\Psi )v)\right)\,,
\eeq
where in the last equality we used the definition of the dual map in the presence of metrics, i.e.~$\kappa(\Psi v,\xi)=g(v,\Psi^*\xi)$ for all $v \in T_xM$ and $\xi \in A_x$.
On the other hand we have
\beqn
0 = \eta_0 (v+I\Psi v, I\xi)=\eta_0 (v, I\xi)+\eta_0 (I\Psi v, I\xi) =\eta_0 (v, \xi)+\eta_0 (\Psi v, \xi)=\eta_0 (v, \xi) +\kappa(\Psi v,\xi)
\eeq
for all $v\in T_xM$, $\xi\in A_x$, where in the third equality we used that  $\Ker (1-I)=TM$.
This gives the following formula for $\eta_0$ in terms of $g$, $\kappa$, and $\Psi$:
\beq\label{metric_eta_0_eqn1}
\eta_0 (v,\xi) &=& -g(v,\Psi^*\xi)\\ \label{metric_eta_0_eqn2}
\eta_0 (\xi,\xi) &=& \kappa(\xi,\xi) \\ \label{metric_eta_0_eqn3}
\eta_0 (v,v) &=& g\left( (1+\Psi^*\Psi)v,v\right)
\eeq
where $v\in T_xM$ and $\xi\in A_x$.

\vskip 2mm\noindent Now let us apply the orthogonality of $I$ one more time. Given that $I$ is an involution, i.e.~$I^2=\Id$,
so that $T\calG|_M=\Ker (1-I)\oplus \Ker (1+I)$,
it is sufficient to require
that the $\pm 1$ eigenspaces of $I$ are orthogonal to one another. From $\Ker (1-I)=TM$, we obtain that one should have $\eta_0 ((1-I)\xi, v)=0$
for any $\xi\in A_x$ and $v\in T_xM$. Using the isomorphism between $A$ and $V(t)$, we identify the anchor map $\rho$ with the restriction of $\md s$ to $V(t)$.
Taking into account that $(1+I)\xi\in T_xM$, one has  $(1+I)\xi = \rho (1+I)\xi$, and since $I\xi\in V_x(s)$, even $ (1+I)\xi =\rho (\xi)$. We therefore obtain the following identity:
\beqn
0=\eta_0 ((1-I)\xi, v)=\eta_0 ((1+I)\xi, v)-2\eta_0 (I\xi, v)=\eta_0 (\rho\xi, v)-2\eta_0 (\xi, v)\,.
\eeq
From Equations (\ref{metric_eta_0_eqn1}) and (\ref{metric_eta_0_eqn3}) we immediately deduce that
\beqn
g((1+\Psi^*\Psi)\rho\xi,v)=-2g(v,\Psi^*\xi)
\eeq
for all $\xi\in A_x$ and $v\in T_xM$. Therefore
\beq\label{eqn:for_Psi}
(1+\Psi^*\Psi)\rho =-2\Psi^*\,.
\eeq
The identity (\ref{eqn:for_Psi}) implies the following equation for $\Psi\Psi^*$
\beqn
\Psi^*\Psi \left(1+\Psi^*\Psi\right)^{-2}=\left(1+\Psi^*\Psi\right)^{-1} - \left(1+\Psi^*\Psi\right)^{-2}=\frac{1}{4}\rho\rho^*\,,
\eeq
or, equivalently,
\beq\label{eq:necessary_inequality_preliminary}
\left(\left(1+\Psi^*\Psi\right)^{-1}-\frac{1}{2}\right)^2=\frac{1}{4}(1-\rho\rho^*)\,.
\eeq
From (\ref{eq:necessary_inequality_preliminary}) we immediately obtain \eqref{smooth} as a necessary condition for the existence
of $\eta_0$ (and thus of a compatible Riemann-Cartan structure on $\calG$). This proves the existence in the second part of Theorem \ref{thm:bi-invariant_metric}. Under the condition \eqref{smooth},
Equation (\ref{eq:necessary_inequality_preliminary}) implies
\beq\label{eq:minus_preliminary}
\left(1+\Psi^*\Psi\right)^{-1}=\frac{1}{2}\left( 1-  \left(1-\rho\rho^*\right)^{\frac{1}{2}} \right)
\eeq
or
\beq\label{eq:plus_preliminary}
\left(1+\Psi^*\Psi\right)^{-1}=\frac{1}{2}\left( 1+ \left(1-\rho\rho^*\right)^{\frac{1}{2}} \right)\,.
\eeq
Given that $0<\left(1+\Psi^*\Psi\right)^{-1}\le 1$, we obtain an additional constraint $\rho\rho^*>0$  in the
first case (\ref{eq:minus_preliminary}) and no additional constraints in the second case (\ref{eq:plus_preliminary}). Therefore,
the equation (\ref{eq:necessary_inequality_preliminary}) possesses at least the solution
\beqn \label{unique}
\Psi^*\Psi =  \rho\rho^*\left(1+ \left(1-\rho\rho^*\right)^{\frac{1}{2}}\right)^{-2} \, ,
\eeq
resulting from equation (\ref{eq:plus_preliminary}), which,  when combined with (\ref{eqn:for_Psi}),  permits to determine $\Psi$: \beq\label{eq:plus_solution}
\Psi = -\rho^*\left(1+ \left(1-\rho\rho^*\right)^{\frac{1}{2}}\right)^{-1}\, .
\eeq
This proves the existence in the second part of Theorem \ref{thm:bi-invariant_metric}.
It remains to prove the uniqueness properties as formulated there.

\vskip 2mm\noindent In regions where $\rho\rho^*>0$, we can get another local solution for $\Psi$ and thus for $\eta$  by using equation (\ref{eq:plus_preliminary}); denoting the previously obtained solution \eqref{unique} by $\Psi_+$ in what follows and this second one by  $\Psi_-$, we have the following two local solutions:
\beq\label{eq:minus_solution}
\Psi_{\pm} = -\rho^*\left(1\pm \left(1-\rho\rho^*\right)^{\frac{1}{2}}\right)^{-1} \, .
\eeq
The condition that $\rho\rho^*$ has no zero eigenvalues, $\rho\rho^*>0$, is equivalent to $A$ being a transitive Lie algebroid. (By definition, $A$ is transitive, if, for all $x \in M$, $\rho_x$ is surjective or, equivalently, $\rho^*_x$ is injective).

\vskip 2mm\noindent From (\ref{eq:minus_preliminary}) and (\ref{eq:plus_preliminary})
  we conclude that if $\left(\Psi_+\right)_x$ and $\left(\Psi_-\right)_x$ exist and coincide at some point $x\in M$, then
   one has $\rho_x\rho_x^*=1$; the proof of the uniqueness will be based upon this fact as well as upon the following Lemma \ref{lem:transitive_extension}.
     \begin{lemma}\label{lem:transitive_extension}
  Let $\calG$ be a Cartan groupoid over $M$, $(A, \rho, g, \kappa)$ be the Lie algebroid of $\calG$, endowed with a positive quadratic structure,
  and $N$ be a path-connected subset of $M$, such that $\rho$ is surjective
  at each point of $N$. Let $\rho_x\rho_x^*=1$ at some point $x\in N$, then
  $\rho\rho^*\equiv 1$ over $N$.
  \end{lemma}
  \noindent\proof(of Lemma \ref{lem:transitive_extension}) Given any point $y\in N$, there exists a path $\gamma$ connecting $y$
         and $x$. Let $\tilde\gamma$ be an $A-$path which lifts $\gamma$, $\mathrm{d} \gamma = \rho\circ\tilde\gamma$, which always exists since
         $A$ is transitive over $N$. By exponentiation of $\tilde\gamma$ we obtain a path on the $t-$fiber over $y$,
         which starts at the identity at $y$ and ends at some element $u\in\calG$ such that $s(u)=x$, $t(u)=y$.
         Taking into account
         that both $\rho$ and $\rho^*$ are invariant under the Cartan groupoid adjoint action and thus
         under $Ad_u$, we immediately deduce  that $\rho_y\rho_y^*=1$.
  $\square$

  \begin{cor}[of Lemma \ref{lem:transitive_extension}]\label{cor:disjoint_union} Let $A$ be as in Lemma \ref{lem:transitive_extension}, $M_1$ be the set of points where $\rho\rho^*$ is equal to $1$ and $M_2$ its complement, $M_2 = M \setminus M_1$. Then $M=M_1 \cup M_2$ is a disjoint union of topological spaces. In particular, if $M$ is connected, either $\rho\rho^*=1$ everywhere or nowhere.
  \end{cor}
  \noindent\proof(of Corollary \ref{cor:disjoint_union}) Suppose $x$ is an arbitrary point of $M_1$, then $\rho_x$ is necessarily surjective; this asserts that there exists a path-connected open neighborhood $U_x$ of $x$ such that $\rho$ is still surjective over $U_x$. By Lemma \ref{lem:transitive_extension} we have $\rho_y\rho_y*= 1$ for all $y\in U_x$, which means that the neighborhood $U_x$ is contained in $M_1$. By the definition of an open set (each point containing an open neighbourhood), $M_1$ is open. On the other hand, $M_1$ is closed as the preimage of a closed set (the unit section of the endomorphism bundle of $TM$) with respect to the continuous function $\rho\rho^*$. Therefore also $M_2$ is open and closed, which implies the statement.  $\square$

\vskip 2mm\noindent Now we are ready to conclude the proof. According to the above lemma and its corollary, provided $M$ is connected, if there is one point $x \in M$ such that $\rho_x\rho_x*= 1$, then $\rho\rho^*=1$ everywhere on $M$ and there is a unique solution for $\eta$ resulting from $\Psi = \Psi_\pm = - \rho^*$. In all other cases, two local solutions $\Psi_+$ and $\Psi_-$ cannot be glued together. Thus, for the second solution $\Psi\ne \Psi_+$ to exist globally on a connected base, one needs $\rho\rho^*>0$ everywhere on $M$, which is equivalent to $A$ being transitive. For this second solution to be different from   $\Psi_+$ we still need to exclude the case where the transitive Lie algebroid $A$ is such that $\rho\rho^*=1$ at one point and thus everywhere. $\blacksquare$

\vskip 2mm\noindent Provided the base $M$ is compact, the norm of $\rho\rho^*$ is necessarily bounded, thus the above-mentioned constraint (\ref{smooth}) can be achieved by proper rescaling of the metric $g$ or $\kappa$ (evidently, $\rho\rho^*<1$ is sufficient for satisfying  (\ref{smooth})). Note that in the gauge theoretic applications, this would be achieved by a sufficiently \emph{small} coupling constant, which is necessary for the application of perturbative methods in any case.

\vskip 2mm\noindent  In the non-compact case, we can bypass the constraint by choosing the initial metric on
$T\calG|_M =TM\oplus V(t)|_M$ as
\beq\label{eq:initial_metric}
\eta_{\sst in}=g\oplus \kappa\,;
\eeq
it will be $Ad_\calG$-invariant, but, in general,  not $I$-invariant. Taking subsequently  the average $\eta_0 =\frac{1}{2}\left(\eta_{\sst in} +I^*\eta_{\sst in} \right)$, the obtained metric remains $Ad_\calG$-invariant, but now is also
invariant under the inversion. However, the obtained Riemann-Cartan structure on $\calG$ will not be compatible with $\kappa$ and $g$ in the sense of Definition \ref{def:compatible_CLR_metric}. We will come back to such ideas once more below.

\begin{rem}\label{rem:direct_metric}Let us assume that the condition (\ref{smooth}) is satisfied, then, as a somewhat longer calculation shows, if we replace $g$ and $\kappa$ in (\ref{eq:initial_metric}) with the also $Ad_\calG$-invariant metrics
$g(R\cdot, \cdot)$ and $\kappa (C\cdot, \cdot)$, respectively,
where $R=2\left(1+\left(1-\rho\rho^*\right)^{\frac{1}{2}}\right)^{-1}$and $C=1-\frac{1}{2}\rho^*R\rho$, we will
get the solution corresponding to (\ref{eq:plus_solution}).
\end{rem}

\begin{example}\label{ex:action_groupoid}
Let $(M,g)$ be a Riemannian manifold, $G$ be a compact Lie group, acting on $M$ by isometries and let $\calG=G\ltimes M$ be the corresponding action Lie groupoid together with the canonical
flat Cartan-Lie structure, chosen such that the leaves of the underlined multiplicative distribution are "constant bisections", i.e.~fibers of the projection $\calG\to G$.
Let $\g$ be the Lie algebra of $G$, endowed with some bi-invariant Riemannian metric, and $\kappa$ be the corresponding flat metric on the action Lie algebroid $A\simeq \g\times M$. By Remark \ref{rem:direct_metric}, we can supply $\calG$ with a compatible Riemann-Cartan structure provided the condition (\ref{smooth}) holds.
\end{example}

\begin{example}\label{ex:action_groupoid_linear}
Take the setting of the last example for the case of $(M,g)$ being an Euclidean vector space $(\mathbb{V}, g_\mathbb{V})$, i.e.~$M=\mathbb{V}$ is a vector space and $g$ then a constant metric tensor on it induced by the bilinear form  $g_\mathbb{V}$. Then the necessary condition for compatibility (in the natural sense of Definition \ref{def:compatible_CLR_metric}), $1-\rho \rho^* \geq 0$, \emph{never} holds true, unless $\rho=0$: with $\rho$ being linear, $\rho \rho^*$ is quadratic and thus unbounded over $M=\mathbb{V}$.
\end{example}

\vskip 2mm\noindent There are several ways to bypass the exclusion of examples as the previous one, as may be useful to know for example in the context of gauge theoretic applications, each of which having different advantages and disadvantages. One of those was already mentioned around equation \eqref{eq:initial_metric}. A likewise option is to equip $T\calG|_M =TM\oplus V(s)|_M$ with $g \oplus \kappa$ and extend it to all of  $\calG$ so as to obtain a bi-invariant metric. In the case of the action groupoid this leads to the natural Riemannian product structure $(G,\hat \kappa) \times (M,g)$, where $\hat \kappa$ is the left-extension of $\kappa$ to the group. Now the unit section is even an isometric embedding, which it cannot be in the case of Definition \ref{def:compatible_CLR_metric}. Again this Riemannian structure on $\calG$ is not inversion invariant and the $t$-map is a Riemannian submersion, but for another metric on the base, namely for $\tilde g = g (\tilde{R}\cdot ,\cdot )$, where $\tilde{R}=\left(1+\rho\rho^*\right)^{-1}$. Coming back to the special case of an action groupoid, there is a second natural way of identifying the total groupoid manifold  as a product, namely so that the $t$-map becomes the projection to the factor $M$. Equipping this splitting of $\calG$ with the product Riemannian structure, one obtains a \emph{different} metric on $\calG$ than the one above;  this second choice corresponds to \eqref{eq:initial_metric} certainly.

\vskip 2mm\noindent In the above two procedures, it was either the $t$-fibers or the $s$-fibers, both at the identity section, that were metrically identified with the fibers of $(A,\kappa)$  and then declared orthogonal to the $TM$-fibers; cf.~also Lemma \ref{lem:adjoint_action}. The respective other fiber then carries a different metric; for example, in the case of the $t$-fiber identification \eqref{eq:initial_metric}, $V(s)|_M$ will be equipped with the metric $\tilde \kappa = \kappa(\tilde{C}\cdot,\cdot)$, where $\tilde{C}=1+\rho^*\rho$; $\kappa \neq \tilde \kappa$ is a direct reflection of the fact that the resulting groupoid is not Riemannian with such a choice.

  \vskip 2mm\noindent There is another option related to the above ones in a kind of democratic way that now \emph{does} produce a Riemannian groupoid. We already mentioned in the proof of Theorem \ref{thm:bi-invariant_metric} that for a Riemannian Cartan groupoid the $\pm 1$-eigenspaces of $I$, $V_\pm \subset T\calG|_M$, must be orthogonal to one another and that $V_+ = TM$. We can now use this last identification to equip $V_+$, a bundle over the identity section $M$, with the fiber metric induced by the  Riemannian manifold $(M,g)$. The projector to $V_-$ is $(1-I)/2$;
using that $(1+I)\xi = \rho(\xi)$ for any $\xi \in V(t)\vert_M$, we find the identification of elements of $V(t)\vert_M$ with vectors from $V_-$ by means $\xi \mapsto (1-\rho/2) \xi$. This isomorphism permits us to equip $V_-$ with a fiber metric as well, pushing forward $\kappa$ on $A$. By requiring that the $V_+$ and $V_-$ are orthogonal, we obtain a metric on $T\calG|_M$  invariant under the involution $I$.
Provided both $g$ on $M$ and $\kappa$ on $A$ are $Ad_\calG-$invariant, we immediately extend the corresponding $Ad_\calG-$invariant fiber metric on $T\calG|_M$ to a Riemann-Cartan structure on $\calG$ by parallel transport.
By construction, the induced metric on the identity section inside $\calG$ is $g$. On the other hand, the $s$- and $t$-fibers do not carry the metric $\kappa$, but $\kappa' = \kappa (C'\cdot ,\cdot)$,
where $C'=1+\frac{1}{4}\rho^*\rho$ and likewise, the metric on the base induced by the Riemannian submersions $\mathrm{d}s$ and $\mathrm{d}t$ is $g' = g(R'\cdot,\cdot)$, where $R'=(1+\frac{1}{4}\rho\rho^*)^{-1}$. While the $Ad_\calG-$invariant metrics $g$ and $\kappa$ are arbitrary, the second couple of $Ad_\calG-$invariant metrics $g'$ and $\kappa'$, which is now compatible with the Riemann-Cartan structure on $\calG$ in the sense of Definition \ref{def:compatible_CLR_metric}, will satisfy the condition
(\ref{smooth}), where the conjugate of $\rho$ is taken with respect to $g'$ and $\kappa'$;  although this condition follows from Theorem \ref{thm:bi-invariant_metric}, one can also easily verify it directly.

\vskip 2mm\noindent So this last ''democratic choice'' has the advantage of not excluding examples such as Example \ref{ex:action_groupoid_linear}, while still staying inside the realm of Riemannian groupoids. However, the price to be paid is that the  usual identification of $V(t)\vert_M$ with $A$  does neither respect the fiber metric nor the one on the base. Now the groupoid is not over $(M,g)$  anymore but over $(M, g')$ (while on the other hand the unit section now is indeed $(M,g)$). Note as an aside also, that with this choice applied to the action groupoid, the total manifold does in general no more split into the product of two Riemannian manifolds anymore.

\vskip 2mm\noindent Starting from a quadratic Lie algebroid there are several options to arrive at a Cartan groupoid with a Riemmanian metric on it, satisfying ''appropriate'' compatibility conditions. If one has particular applications in mind, one of the alternative choices just mentioned above might be possibly more adequate. Note that for the above-mentioned alternatives there are no obstructions whatsoever, and the construction is in each case rather straightforward. On the other hand, we found in this paper with Theorem \ref{thm:bi-invariant_metric} the optimal statement that one can make given the natural setting of a Riemannian Cartan-Lie groupoid as defined in Definition \ref{deff:Cartan-Lie_Riemannian_groupoid} and the natural identification of a quadratic Lie algebroid corresponding to it as in Definition \ref{def:compatible_CLR_metric}.

\begin{rem}\label{rem:Riemannian_general} By construction, a bi-invariant metric $\eta$ on a Cartan groupoid $\calG$, invariant under the inversion map, is a particular case of a
general Riemannian groupoid (see \cite{GGHR1989}) or a groupoid $1$-metric in terms of \cite{delHoyo-Fernandes15}. A natural question is whether
$\eta$ satisfies a stronger property of the extension to a $2$-metric, so that a Riemann-Cartan groupoid is a Riemannian groupoid in the sense of \cite{delHoyo-Fernandes15} (see also \cite{delHoyo-Fernandes16})? If a certain sufficient condition which we will specify below is fulfilled, a ($\calG-$invariant) $2$-metric extension of a Riemann-Cartan structure can be derived from the proof of Theorem 4.3.4 in \cite{delHoyo-Fernandes15}. The authors of \cite{delHoyo-Fernandes15} use a cotangent average with respect to a Haar density on a proper Hausdorff groupoid to make all face maps of the classifying space of $\calG$ into Riemannian submersions. On the other hand, if we start with any bi-invariant metric on a Cartan groupoid, the cotangent averaging will be not required; we just skip the first step of their proof and follow the rest of it. Thus we obtain a (fully extended) metric on the nerve $N^\bullet\calG$, such that
$N^k\calG$ is invariant under the action of the symmetric group $S_{k+1}$ for each $k\ge 1$. In particular, we get a metric on $\calG =N^1\calG$, which is
bi-invariant and inversion invariant, simultaneously; it remains to check when a Riemann-Cartan structure $\eta$ can be obtained as a metric on $N^1\calG$ in this way,
i.e.~as a part of a fully extended metric on $N^\bullet\calG$. Assume that the restriction of $\eta$
onto (above-mentioned) $V_+$ and $V_-$ coincides with $g$ and $\kappa$, respectively. Provided $V_+$ and $V_-$ are identified  with $TM$ and $A$, the sufficient condition for $(\calG,\eta)$
to carry a $\calG-$invariant fully extended metric---and thus, in particular, a 2-metric---is $\rho\rho^*<4$, where the conjugation is taken with respect to $g$ and $\kappa$. Note that in this setting there were no conditions on $g$ and $\kappa$ to be satisfied for the Riemann-Cartan structure $\eta$ to exist, cf. the considerations preceding this remark; so this condition will not be a consequence of the condition for the integration found in Theorem  \ref{thm:bi-invariant_metric}, but an additional, stronger one, also when formulated in terms of the original data of that theorem.
\end{rem}

\end{document}